\documentclass[a4paper,12pt,reqno]{amsart}
\usepackage{amsfonts}

\usepackage[T1]{fontenc} 
\usepackage[latin1]{inputenc}

\pdfoutput=1

\usepackage{tikz}
\usepackage{latexsym}
\usepackage[colorlinks=true, pdfstartview=FitV, linkcolor=blue, citecolor=blue, urlcolor=blue,pagebackref=false]{hyperref}

\usepackage{amsmath, amsthm, amssymb, color, graphicx, hyperref, mathrsfs}

\usepackage[margin=1.3in,marginparwidth=1.5cm, marginparsep=0.5cm]{geometry}

\usepackage{chngcntr} 

\usepackage{booktabs}
\usepackage{multirow}
\usepackage{siunitx}

\usepackage{microtype}

\usepackage{mathrsfs}

\usepackage{color}

\usepackage{marginnote}
\usepackage{comment}

\definecolor{darkgreen}{rgb}{0,0.5,0}
\definecolor{darkred}{rgb}{0.7,0,0}
\definecolor{darkblue}{rgb}{0,0,0.7}
\theoremstyle{plain}
\newtheorem{theorem}{Theorem}[section]
\newtheorem{corollary}[theorem]{Corollary}
\newtheorem{lemma}[theorem]{Lemma}
\newtheorem{proposition}[theorem]{Proposition}

\newtheorem{definition}[theorem]{Definition}

\newtheorem{conjecture}[theorem]{Conjecture}
\newtheorem{question}[theorem]{Question}

\theoremstyle{remark}
\newtheorem{remark}[theorem]{Remark}

\numberwithin{equation}{section}
\numberwithin{table}{section}

\newcommand{\N}{\mathbb{N}}
\newcommand{\R}{\mathbb{R}}

\newcommand{\Z}{\mathbb{Z}}
\newcommand{\Q}{\mathbb{Q}}

\newcommand{\E}{\mathbb{E}}
\newcommand{\T}{\mathbb{T}}

\renewcommand{\Q}{\mathbb{Q}}

\renewcommand{\Z}{\mathbb{Z}}
\renewcommand{\N}{\mathbb{N}}

\def\T{\mathbb{T}}

\let\qed=\QED

\def\calG{\mathcal{G}}

\def\calT{\mathcal{T}}

\def\ceil#1{\lceil{#1}\rceil}
\def\floor#1{\lfloor{#1}\rfloor}

\def\P{\mathbb{P}} 
\def\E{\mathbb{E}} 
\def\md{\mid}

\def \eps {\epsilon}

\def\Bb#1#2{{\def\md{\bigm| }#1\bigl[#2\bigr]}}

\def\Pb{\Bb\P}
\def\Eb{\Bb\E}

\def\FK#1#2#3{{\def\md{\bigm| } \P_{#1}^{\,#2}  \bigl[  #3 \bigr]}}
\def\EFK#1#2#3{{\def\md{\bigm| } \E_{#1}^{\,#2}  \bigl[  #3 \bigr]}}

\def \p {{\partial}}

\renewcommand{\subset}{\subseteq}

\newcommand{\F}{\mathscr{F}}

\renewcommand{\hat}{\widehat}

\def\<#1{\langle #1\rangle}

\def\nn{\nonumber}
\def\bi{\begin{itemize}}  
\def\ei{\end{itemize}}
\def\bnum{\begin{enumerate}} 
\def\enum{\end{enumerate}}
\def\ni{\noindent}
\def\bf{\bfseries}

\usepackage{mhequ} 

\colorlet{symbols}{blue!90!black}
\colorlet{testcolor}{green!60!black}

\def\1{\mathbf{{1}}}

\def\g2{\frac {\gamma^2} 2}

\def\action{\curvearrowright}

\title[Diffuse-extensive-amenability and the IET group]
{Inverted orbits of exclusion processes, diffuse-extensive-amenability and (non-?)amenability of the interval exchanges} 

\author{Christophe Garban}

\address[Christophe Garban]
{Univ Lyon, Université Claude Bernard Lyon 1, CNRS UMR 5208, Institut Camille Jordan, 69622 Villeurbanne, France}
\email{garban@math.univ-lyon1.fr}
\begin{document}

\maketitle


\begin{abstract}
The recent breakthrough works \cite{Kate,Mik1,Mik2} which established the amenability for new classes of groups,  lead to the following question: is the action $W(\Z^d) \action \Z^d$ extensively amenable? (Where $W(\Z^d)$ is the {\em wobbling group} of permutations $\sigma:\Z^d \to \Z^d$ with bounded range). 
This is equivalent to asking whether the action 
$(\Z/2\Z)^{(\Z^d)} \rtimes W(\Z^d) \action  (\Z/2\Z)^{(\Z^d)}$ is amenable. 
The $d=1$ and $d=2$ and have been settled respectively in \cite{Kate,Mik1}. By \cite{Mik2}, a positive answer to this question 
 would imply the amenability of the IET group. In this work, we give a partial answer to this question by introducing a natural strengthening of the notion of {\em extensive-amenability} which we call {\em diffuse-extensive-amenability}. 

Our main result is that for any bounded degree graph $X$, 
the action $W(X)\action X$ is diffuse-extensively amenable if and only if $X$ is recurrent.
Our proof is based on the construction of suitable stochastic processes $(\tau_t)_{t\geq 0}$ on $W(X)\, <\, \mathfrak{S}(X)$ whose {\em inverted orbits}
$$
\bar O_t(x_0) = 
\{x\in X, \exists s\leq t,\,  \tau_s(x)=x_0\} = \bigcup_{0\leq s \leq t} \tau_s^{-1}(\{x_0\})
$$
are exponentially unlikely to be sub-linear when $X$ is transient. 
This result leads us to conjecture that the action $W(\Z^d) \action \Z^d$ is not extensively amenable when $d\geq 3$ and that a different route towards the (non-?)amenability of the IET group may be needed. 
\end{abstract}

\section{Introduction}

\subsection{IET group, wobbling group $W(\Z^d)$ and criteria of amenability}

The IET group (=group of Interval Exchanges Transformations) is the group of càdlàg  piecewise-translation bijections $g$ from $\T=\R/\Z \to \R /\Z$ s.t. the set 
\[
t(g):= \{g(x)-x \text{ mod } 1, x\in \T\} \subset \T
\]
is finite. Interval Exchanges Transformations have been  used extensively in the realm of dynamical systems, for example in the analysis of  polygonal billards with rational angles. See the survey by Viana \cite{Viana}.  

Katok asked long ago whether the free group $\mathbb{F}_2$ can be embedded as a sub-group in IET. Even though the answer to this question remains elusive,  some progress has been done over the last decade: for example it has been proved in \cite{Dahmani} that for {\em generic} pairs $g_1,g_2 \in \mathrm{IET}$, the sub-group $\<{g_1,g_2}<\mathrm{IET}$ is not free. Another natural open problem in this area is to ask the following question (see \cite{Cornulier}): 

\begin{question}
Is the IET group amenable ?
\end{question}

Note that one does not endow the IET group with any topology here. Therefore amenability is equivalent to the fact that all finitely generated sub-groups $H=\<{g_1,\ldots,g_k} < \mathrm{IET}$  are amenable. 

We will use throughout the following key fact: for any  $\{g_1,\ldots,g_k\} \subset \mathrm{IET}$, the sub-group $H=\<{g_1,\ldots, g_k}$ can be embedded in the group of permutations of an Euclidean lattice $\Z^{\bar d} \times \Z/m \Z$. The parameters $\{\bar d ,m\}$ are uniquely determined by $g_1,\ldots, g_k$ and $\mathrm{rk}_{\Q}(g_1,\ldots,g_k):= \bar d$ is called the \textbf{rational rank} of the sub-group $H$. 
Let us briefly outline how it works. 

\bi
\item First, the parameters $\{\bar d,m\}$ are obtained as follows: consider the sub-group $\Lambda< \T=\R/\Z$ defined by 
\[
\Lambda = \Lambda(H):= \<{t(g_1) \cup t(g_2) \cup \ldots \cup t(g_k)}  \subset \T\,.
\]
Note that $\Lambda$ is defined equivalently as generated by $\bigcup_{g\in H} t(g)$ and thus does not depend on the choice of generators of $H$. 
Abelian sub-groups of the circle  are of the form $\Lambda \simeq \Z^{\bar d} \times \Z/m\Z$ with $m\geq 1$ (if $m=1$, it is just $\Z^{\bar d}$). Equivalently, one can find $\alpha_1,\ldots, \alpha_{\bar d} \in \R$ which are linearly independent over $\Q$ such that 
the map $\Phi(n_1,\ldots,n_{\bar d},b):= \sum n_i \alpha_i +  b \frac 1 {m} \text{ (mod $1$)}$ is an isomorphism from $\Z^{\bar d} \times \Z/ m \Z$ to $\Lambda \subset \T$. 

\item One may now associate explicit permutations $(\sigma_i)_{1\leq i \leq k}$ to each interval exchanges $g_1,\ldots,g_k$ as follows:  for each $1\leq i \leq k$, define
\begin{align*}\label{}
\sigma_i  \,:   \, \left\{
\begin{array}{rcl}
 \Z^{\bar d}  \times    \Z/ m\Z & \to & \Z^{\bar d}   \times    \Z/ m \Z \\
 (\vec{n}, b)& \mapsto &  \Phi^{-1}\circ g_i \circ \Phi (\vec n,  b)
\end{array}
\right.
\end{align*}
\ei
See for example \cite{Mik2}. In order to avoid carrying the finite group $\Z/m\Z$ throughout, note that these permutations can always be viewed as permutations on a (possibly larger) Euclidean lattice $\Z^d$. This is  clear if $m=1$ (by taking $d=\bar d$) and if $m\geq 2$, by identifying $\Z/ m \Z$ with $\{0,\ldots, m-1\} \subset \Z$, we may view $\{\sigma_i\}_{1\leq i \leq k}$ as permutations on a
{\em ``slab graph''} $\Z^{\bar d} \times \{0,\ldots, m-1\}\subset \Z^{d=\bar d +1}$. 


It is easy to check that not only $\sigma_1,\ldots,\sigma_k$ belong to the symmetric group $\mathfrak{S}(\Z^d)$\footnote{Throughout this paper, $\mathfrak{S}(\Z^d)$ is the space of \textbf{all} permutations $\sigma : \Z^d \to \Z^d$, in particular they typically do not have a finite support.} but they also belong to the so-called {\em wobbling group} $W(\Z^d)$ of $\Z^d$ defined for general graphs as follows: 
%
%
%
%

\begin{definition}[See  \cite{wobbling}]\label{d.wobbling}
Let $G=(X,E)$ be a locally-finite connected graph. The \textbf{wobbling group} $W(X)$ is defined as:
\begin{align}
W(X):= \Big\{ 
\sigma \in \mathfrak{S}(X), \text{ s.t. }\mathrm{range}(\sigma):=\sup_{x\in X} \{ d_G(x,\sigma(x)) \} <\infty  \Big\}\,,
\end{align}
where $\mathfrak{S}(X)$ is the group of all permutations $\sigma : X \to X$.

We will equip $W(X)$ (resp. $\mathfrak{S}(X)$) with the topology $\calT$ of pointwise convergence where $\sigma_n \overset{\calT}\longrightarrow \sigma \in W(X)$ (resp.  $\mathfrak{S}(X)$) if and only if for any $x\in X$, $\sigma_n(x) \to \sigma(x)$.
This topology is metrizable, for example using the following distance $d_{x_0}$ defined on $W(X)$ (resp. $\mathfrak{S}(X)$) for any root-vertex $x_0\in X$: let $d_{x_0}(\sigma, \tilde \sigma)=2^{-n}$ if $\sigma_{|B(x_0,n)}=\tilde \sigma_{|B(x_0,n)}$ and $\sigma^{-1}_{|B(x_0,n)}=\tilde \sigma^{-1}_{|B(x_0,n)}$ and at least one of these two equalities fails for the ball $B(x_0,n+1)$. It is easy to check that equipped with this distance, the metric space $(\mathfrak{S}(X),d_{x_0})$ is a Polish space. 

Let us introduce a family of compact subsets $W_r(X)$ of $W(X)$ indexed by $r\in \N_+$ which will be used throughout. For any $r\geq 1$, let 
\begin{align}\label{}
W_r(X):= \{ \sigma \in W(X), \text{ s.t. }\mathrm{range}(\sigma) \leq r \}\,.
\end{align}
(Note that $W_r(X)$ is not a subgroup of $W(X)$).
\end{definition}

\begin{remark}\label{r.wobbling} $ $ 
\bnum
\item If one prefers to stick to permutations on the abelian group $\Lambda \simeq  \Z^{\bar d} \times  \Z/m \Z$, the permutations $\sigma_i$ would leave in the wobbling group $W(X=\Z^{\bar d} \times \Z/m \Z)$ for the natural graph structure on $X$. 
\item The wobbling groups $W(\Z^d)$ have been studied recently  for example in \cite{wobbling}.
\item 
There exist finitely generated sub-groups of $W(\Z^d)$ that cannot be obtained from the group IET through the above procedure.
See Remark \ref{r.criterion} below. 
\enum
\end{remark}

The above discussion provides us with an explicit embedding of
$H=\<{g_1,\ldots,g_k}<\mathrm{IET}$ into $\<{\sigma_1,\ldots, \sigma_k}< W(\Z^d)$. As such the amenability of IET is equivalent to the amenability of every  finitely generated subgroups of $\{W(\Z^d), d\geq 1 \}$.


In what follows, we shall fix some $\{g_1,\ldots, g_k\}$ in IET and we will investigate criteria of amenability for its associated group of permutations $G:=\<{\sigma_1,\ldots, \sigma_k}< W(\Z^d)$. 
In order to list some criteria for the amenability of $G < W(\Z^d)$, it is useful to introduce the Cayley graph $\Gamma$ on $G$ for the symmetric generating set $S=\{ \sigma_1^{\pm 1}, \ldots, \sigma_k^{\pm 1}\}$.

\smallskip

\ni
\textbf{List of equivalent criteria of amenability}. We list below some criteria of amenability in our special case of $G=\<{\sigma_1,\ldots, \sigma_k}<W(\Z^d)$ but the equivalence of these  criteria holds much more generally for any finitely generated groups. (See for example \cite{Paterson} for a complete list of criteria).
\bnum
\item \textbf{$G$-invariant mean.}
There exists a finite and finitely additive measure $m$ on $G$ which is $G$-invariant.
\item \textbf{Fölner sequence.}  There exists a sequence of finite subsets $F_n \nearrow G$, s.t. $|\p F_n|/ |F_n| \to 0$ (where the boundary of a set $J$ is, say, the edge boundary of $J$ in the Cayley graph $\Gamma$). 
\item \textbf{Kesten (version 1).}  There exists a {\em symmetric} probability measure $\mu$ supported on a \underline{finite} generating set of $G$ which is such that if $\{s_i\}_{i\geq 1}$ are i.i.d random variables in $W(\Z^d)$ with $s_i \sim \mu$, then for any $\eps>0$, the random walk $\tau_n:=s_n \circ s_{n-1} \circ \ldots s_1$ returns to the identity (i.e. $\mathrm{Id} : \Z^d \to \Z^d$) at even times $2n$ with probability $p_{2n}$  larger than $e^{-\eps n}$ when $n$ is large enough. 
\item \textbf{Kesten (version 2).} The above estimate on the return probability of (symmetric) random walks holds for any symmetric probability measure $\mu$ on a \underline{finite} subset of $G$. 
\enum

In particular, the amenability of IET is equivalent to proving that for any fixed  $\{g_1,\ldots,g_k\}\subset \mathrm{IET}$, the random walk on permutations in $W(\Z^d)<\mathfrak{S}(\Z^d)$ induced by the uniform measure $\mu$ on $S=\{ \sigma_i^{\pm 1}\}_{1\leq i \leq k}$ returns to $\mathrm{Id}_{\Z^d\to\Z^d}$ with probability larger than $e^{-\eps n}$ for any $\eps>0$ and $n$ large enough.

\subsection{A new criterion for amenability}
The works \cite{Kate,Mik1,Mik2} introduced a striking method which lead to a set of new useful criteria to establish the amenability a large class of groups. In order to introduce their criterion in the case of IET, we need the following notion of \textbf{inverted orbit.}

\begin{definition}\label{}
Let $(\tau_n)_{n\geq 0}$ be a sequence of elements in $\mathfrak{S}(\Z^d)$ (for example a random walk with values in $\mathfrak{S}(\Z^d)$). For each $n\geq 0$, the \textbf{inverted orbit} of the origin is the subset $O_n \subset \Z^d$ defined as 
\begin{align*}\label{}
O_n := \{x\in \Z^d, \exists j \in \{0,1,\ldots, n\},  \tau_j(x)=0\} = \bigcup_{0\leq j \leq n} \tau_j^{-1}(\{0\})\,.
\end{align*}
\end{definition}

This new criterion from  \cite{Kate,Mik1,Mik2} can be stated as follows for IET (we refer to \cite{Mik2} for other equivalent formulations):
\begin{theorem}[\textbf{``Inverted orbit criterion''} \cite{Mik2}]\label{th.Mik} $ $

\ni
Fix any $\{g_1,\ldots,g_k\} \subset \mathrm{IET}$. Let $\{\sigma_1,\ldots, \sigma_k\}$ be their associated permutations in $W(\Z^d)$, with $d\geq \mathrm{rk}(g_1,\ldots,g_k)$.
 
The subgroup $H=\<{g_1,\ldots,g_k}<\mathrm{IET}$ is amenable if and only if the inverted orbits $\{O_n\}_{n \geq 0}$ of the random walk on $W(\Z^d)$ defined as 
\[
\tau_n := s_n \circ s_{n-1} \circ \ldots s_1\,,\,\,\,   \text{  where $s_i$  are  i.i.d  and } \sim \mu=\sum_{i=1}^k \frac 1 {2k}(\delta_{\sigma_i} + \delta_{\sigma_i^{-1}})
\]
satisfy any of the following equivalent conditions:
\bi
\item[i)] For any $\eps>0$, $\Pb{|O_n| \leq \eps n} \geq e^{-\eps n}$ for $n$ large enough. 
\item[ii)] $\lim_{n\to \infty}\frac 1 n \log \Eb{(\frac 1 2)^{|O_n|}} =0$. 
\ei
\end{theorem}

Note that this new criterion is both very surprising and powerful. Indeed, it is very easy to check that the classical amenability criteria listed above imply this one, but the converse is significantly harder (\cite{Mik1,Mik2}). To highlight this, assume that the criterion \textbf{Kesten (version 2)} holds and let us show that it easily implies the above criterion. Fix $\eps>0$, if $h$ is large enough then the above symmetric random walk $\tau_n = s_n \circ s_{n-1} \circ \ldots s_1$ returns to the identity map at time $n=h$ with probability larger than $e^{-\eps h}$.
Now, let us use the fact that each of the $\sigma_i$ belong to $W(\Z^d)$. If $M:=\max_{1\leq i \leq k} \mathrm{range}(\sigma_i)$, in $h$ steps, sites travel at distance at most $h M$. 
This implies in particular that $O_h \subset B(0, h M)$. The key property to notice here is that on the event $\{\tau_h = \mathrm{Id}_{\Z^d \to \Z^d}\}$ one can repeat the same argument and obtain that now, $O_{2h} \subset B(0, h\dot M)$. We can iterate this bound for any $N=m\cdot h$, with $m\geq 1$: the probability that $\tau_N$ returns to the identity at every times $h, 2h, \ldots, N=(m-1)h$ is larger than $e(^{-\eps h})^m = e^{-\eps N}$. This implies that for any such $N$, one has $\Pb{O_N \subset B(0, h M)}  \geq e^{-\eps N}$. One can thus find a constant $C=C(d,M)<\infty$ s.t. $\Pb{|O_N| \leq C h^d}  \geq e^{-\eps N}$ It is not hard to conclude that item $i)$ holds.  We let the reader appreciate the fact that the reverse direction \textbf{``Inverted orbit criterion'' $\Rightarrow$ Kesten (version 2)} is much less clear. See \cite{Mik1,Mik2}. 

\begin{remark}\label{r.criterion} $ $
\bnum
\item One may wonder whether this amenability criterion holds for any arbitrary subgroup $\<{\sigma_1,\ldots,\sigma_k} < W(\Z^d)$, but this is known not to be the case. Indeed it is shown in \cite{Van} that  $\mathbb{F}_2< W(\Z^d)$ for any $d\geq 1$. If $d=2$, and $\sigma_1,\sigma_2 \in W(\Z^2)$ are s.t. $\<{\sigma_1,\sigma_2} \simeq \mathbb{F}_2$, it follows from \cite{Mik1} that the random walk induced by the uniform measure on $\{\sigma_1^{\pm}, \sigma_2^{\pm}\}$ does satisfy items $i)$ and $ii)$. Yet $\<{\sigma_1,\sigma_2}$ is not amenable. 
This shows that it is important somehow for this {\em inverted orbit criterion} to hold that the set $\{\sigma_1,\ldots, \sigma_k\}$ is produced through the above procedure out of some $\{g_1,\ldots, g_k\}\subset \mathrm{IET}$. 
\item It is not hard to check that items $i)$ and $ii)$ in Theorem \ref{th.Mik} are equivalent. See \cite{Mik2}. The reason why we stated this equivalence is that $i)$ is somewhat reminiscent of Kesten's criterion while $ii)$ has a nice interpretation in terms of a certain \textbf{lamplighter random walk} on $\Z^d$. The usual lamplighter RW on $\{0,1\} \wr \Z^d$ corresponds to a random walker on $\Z^d$ which switches on or off lamps along his way on $\Z^d$. It corresponds to the semi-direct product $(\Z/2\Z)^{(\Z^d)} \rtimes \Z^d$ and this group is well known to be \textbf{amenable}.  
The present situation corresponds to another kind of lamplighter random walk on $\Z^d$: at each step, either all lamps switched on in $\Z^d$ simultaneously move to a new site according to a random permutation uniformly chosen among $S=\{\sigma_i^{\pm 1}\}_{1\leq i \leq k}$, or a lamp is added/removed at the origin $0_{\Z^d}$. See \cite{Mik2} for details. For an appropriate choice of finite generating set of $(\Z/2\Z)^{(\Z^d)}\rtimes \<{\sigma_1,\ldots, \sigma_k}$, the probability of return to the identity $(\bar 0)\times \mathrm{Id}_{\Z^d \to \Z^d}$ is exactly $\Eb{(1/2)^{|O_n|}}$. As opposed to the classical lamplighter groups whose amenability is easy to show, the amenability of IET is equivalent to the amenability of these different types of lamplighter groups $(\Z/2\Z)^{(\Z^d)}\rtimes \<{\sigma_1,\ldots, \sigma_k}$.  Such {\em permutation wreath products} were studied for example in \cite{Bartho1,Bartho2,Amir, Virag,Saloff}.
\enum
\end{remark}

\subsection{Random walks in conductances}
Before going further, let us briefly introduce in this subsection a classical type of random walk on $\Z^d$ which will be used throughout this work. 
\begin{definition}[Random walk among conductances]\label{d.conductance}
Let $C: \Z^d \times \Z^d \to \R_+$ be a symmetric kernel satisfying for any $x\in \Z^d$,  
\[
C(x):= \sum_y C(x,y)\in (0,\infty)
\]
Such an operator induces a random walk $\{X_n\}_{n\geq 0}$ on $\Z^d$ defined as the Markov chain on $\Z^d$ satisfying for any $x,y$ and $n\geq 0$, 
\[
\Pb{X_{n+1}=y \md X_n=x}:= \frac{C(x,y)} {C(x)} 
\]
The process $\{X_n\}_{n\geq 0}$ is called the \textbf{random walk in conductances $\{C(x,y)\}_{x,y}$}. (Notice we may well have $C(x,x)>0$ in which case the process has a positive probability to stay at $x$). 
\end{definition}

The main reason why this particular type of RW will be relevant to us is because of the following fact. For any fixed choice of $\{\sigma_1,\ldots, \sigma_k\}\subset \mathfrak{S}(\Z^d)$, if $(\tau_n)_{n\geq 0}$ is the random walk on $\mathfrak{S}(\Z^d)$ defined as $\tau_n =s_n \circ \ldots s_1$ where $s_i$ are uniform in $S=\{\sigma_i^{\pm 1}\}_{1\leq i \leq k}$, then for any initial site $x\in \Z^d$, the stochastic process $n\mapsto \tau_n(\{x\})$ is a random walk in conductances given by 
\begin{align}\label{e.C}
C(x,y):= |\{ 1\leq i \leq k \,: \, \sigma_i(x) = y\}| +|\{  1\leq i \leq k \,: \sigma_i^{-1}(x) = y\}|\,,
\end{align}
for any $x,y\in \Z^d$.
Clearly $C(x,y)=C(y,x)$ and if all $\{\sigma_i\}_{1\leq i \leq k}$ are in the wobbling group $W(\Z^d)$ then $C$ has \textbf{bounded-range} i.e. $C(x,y)=0$ if $\|x-y\|_2 > h$ for some finite range $h$.

For later use, we collect below some classical properties of random walks in conductances.
\begin{proposition}[Classical facts
, see for example \cite{LyonsPeres} or \cite{Kumagai}]\label{pr.conductance} $ $

Let $(X_n)_{n\geq 0}$ be a RW among the conductances $\{C(x,y)\}_{x,y\in \Z^d}$. 
\bi
\item[a)] $(X_n)_{n\geq 0}$ is reversible for the measure $\mu$ on $\Z^d$ given by $\mu(x):=C(x)$ for any $x\in \Z^d$. \footnote{N.B. Reciprocally, any Markov chain on $\Z^d$ which is reversible for a non-degenerate measure $\mu$ is a random walk in conductances.}
\item[b)] If $C$ has bounded-range and if $d\in \{1,2\}$, then $(X_n)$ is \textbf{recurrent}. 
\item[c)] Same conclusion if $C$ has bounded-range on a slab-graph $\Z^{\bar d} \times \{0,\ldots, m-1\}$ with $\bar d \in \{1,2\}$. 
\item[d)] If $C$ has bounded range on $\Z^d$, then there exists $c>0$ s.t. for all $x\in \Z^d$ and all $k\geq 1$, $\FK{}{x}{X_{2k}=x} \geq c k^{-d/2}$. 
\item[e)] If there exists $\delta>0$ s.t. for any neighbouring sites $x,y\in \Z^d$ (i.e. such that $\|x-y\|_2=1$) $C(x,y) \geq \delta$, then if $d\geq 3$, $(X_n)$ is \textbf{transient}.
\ei
\end{proposition}

\subsection{Main question addressed in this work}

Recall the above ``Inverted orbit criterion'' from \cite{Mik1,Mik2} is much sharper (at least in the case of IET) than classical amenability criteria.  One may then argue/hope  that all specificities of the IET group have already been used in the proof in \cite{Mik1,Mik2} that 
\[
\text{\textbf{``Inverted orbit criterion'' $\Rightarrow$ amenability of IET}}
\]
and that only ``generic'' arguments are needed to conclude for the amenability of IET. For example if the answer to the question below turns out to be positive, it would immediately imply that IET is amenable. 
\begin{question}[See question 1.2 from \cite{wobbling}]\label{q.WZd} $ $

For any $d\geq 1$, prove that any of the following equivalent properties holds:  
\bnum
\item The action $W(\Z^d) \action \Z^d$ is extensively amenable.
\item The action $(\Z/2\Z)^{(\Z^d)} \rtimes W(\Z^d) \action (\Z/2\Z)^{(\Z^d)}$ is amenable.
\item For any finite $\{\sigma_1,\ldots, \sigma_k\} \subset W(\Z^d)$, the \textbf{inverted orbits} $\{O_n\}_{n\geq 0}$ of the random walk $(\tau_n)_{n\geq 0}$ on $W(\Z^d)$ generated by the uniform measure on $S=\{\sigma_i^{\pm 1}\}_{1\leq i \leq k}$ satisfy items $i)$ or $ii)$ in Theorem \ref{th.Mik}.
\enum
\end{question}

To our knowledge, up to now, all natural examples of stochastic processes $(\tau_n)_{n\geq 0}$ in $W(\Z^d)$ whose inverted orbits have been successfully analyzed give support to a positive answer to the above question.
%
%
Let us list the situations that have been analyzed so far.
\bnum
\item The main technical tool in the work \cite{Kate} is precisely to answer this question when $d=1$: it is proved that the action $W(\Z) \action \Z$ is extensively amenable. This property is then used in \cite{Kate} to prove the amenability of the {\em topological full group} of a Cantor system. (This gave the first example of a finitely generated simple group which is amenable). 

\item This analysis is extended in \cite{Mik1} where it is shown that $W(\Z^2) \action \Z^2$ is extensively amenable as well. More generally, it is shown that $W(X) \action X$ is extensively amenable if $X$ is locally-finite and recurrent. 
Using this property, the amenability of a large class of groups has is established in \cite{Mik1}.

\item One particular instance of this situation (analyzed in \cite{Mik2}) is the case of random walks on permutations induced by interval exchanges $\{g_1,\ldots,g_k\}$ whose rational rank statisfies $\mathrm{rk}_{\Q}(g_1,\ldots,g_k)\in \{1,2\}$. These lead to random walks $(\tau_n)_{n\geq 0}$ in $W(\Z^{\bar d} \times \{0,\ldots,m-1\})$ with $\bar d \in \{1,2\}$. For any initial state $x$, the processes $n\mapsto \tau_n(\{x\})$ are random walks among bounded-range conductances on $W(\Z^{\bar d} \times \{0,\ldots,m-1\})$ and are therefore recurrent by item $c)$ above. Using this recurrence property, it is shown in \cite{Mik2} that the inverted orbits of such processes satisfy  $\Eb{|O_n|}=o(n)$ which in particular implies $i)$ in Theorem \ref{th.Mik}. Consequently by Theorem \ref{th.Mik}, these subgroups $H=\<{g_1,\ldots,g_k}$ of low rational rank are amenable!
\item Another relevant example is given by the class of \textbf{rigid permutations} on $\Z^d$, given by $\sigma_i^{\mathrm{rigid}}(x):=x+e_i$ where $(e_i)_{1\leq i \leq d}$ is the canonical basis of $\Z^d$. It is not hard to show\footnote{by an argument very similar to  \textbf{ Kesten, version 2  $\Rightarrow$ ``Inverted orbit criterion''}}
that items $i)$ and $ii)$ also hold for this specific choice of rigid permutations $\{\sigma_1^{\mathrm{rigid}},\ldots, \sigma_d^{\mathrm{rigid}}\}$.
\item Our last example deals with asking the analogous question for \textbf{forward orbits} rather than \textbf{inverted orbits}. If $(\tau_n)_{n\geq 0}$ is, say, the random walk on $W(\Z^d)$ defined as $\tau_n =s_n \circ \ldots s_1$ where $s_i$ are uniform in $S=\{\sigma_i^{\pm 1}\}_{1\leq i \leq k}$ for some fixed choice of $\{\sigma_1,\ldots, \sigma_k\}\subset W(\Z^d)$, then as we noticed above the \textbf{forward orbit} of the origin:\begin{align*}\label{}
\vec{O_n}:= \bigcup_{0\leq j \leq n} \tau_j(\{ 0\})\,, 
\end{align*}
is the nothing but the \textbf{range} of a random walk on $\Z^d$ with bounded-range conductances given by ~\eqref{e.C}. In particular we see here that the geometry of \textbf{forward orbits} is much easier to apprehend than the geometry of
\textbf{inverted orbits}\footnote{In fact, \textbf{inverted orbits} naturally show up for the study of the lamplighter RW governed by $(\Z/2\Z)^{(\Z^d)}\rtimes \<{\sigma_1,\ldots, \sigma_k}$, while the easier \textbf{forward orbits} (=range of RW) naturally show up for the study of the classical lamplighter RW corresponding to $(\Z/2\Z)^{(\Z^d)}\rtimes \Z^d$}. (See also the discussion in \cite{Virag}). 
Using classical heat-kernel bounds for such  conductance-random walks (item $d)$ above), we obtain that the probability of the event  $\{\tau_{2k}(0)=0\}$ is greater than $c\, k^{-d/2}$, where $c=c(\{\sigma_1,\ldots, \sigma_k\})>0$. This plus the fact that $\max_i \mathrm{range}(\sigma_i) < \infty$ easily implies that for any dimension $d\geq 1$, one has for any $\eps>0$, 
\begin{align*}\label{}
\Pb{|\vec{O_n}| \leq \eps n} \geq e^{-\eps n}\,, \text{  for $n$ large enough.}
\end{align*}
\enum

\begin{remark}\label{}
Note at this stage that the last two items do not suggest that the dichotomy recurrence/transience should play any particular role in deciding whether the action $W(\Z^d) \action \Z^d$ will be extensively amenable or not.
\end{remark}

%

Our goal in this paper is to go beyond these examples and to analyze what happens for more ``generic'' examples of diffusions with values in $W(\Z^d)$.
Indeed, as the answer to Question \ref{q.WZd} remains open when $d\geq 3$, we will consider the exact same question but for a broader class of random walks than the random walks defined in Question \ref{q.WZd}, item (3). 
To start with we will not even restrict ourselves to diffusions in the {\em wobbling group} $W(\Z^d)$ but we will start by investigating Question \ref{q.WZd} for random walks on the larger symmetric group $\mathfrak{S}(\Z^d)$.   
Without requiring some natural axioms (see the classes A),B),C),D) in Definition \ref{d.classes} below), it is not hard to come up with examples of random walks on $\mathfrak{S}(\Z^d)$ with i.i.d increments which do not satisfy items $i)$ or $ii)$ from Theorem \ref{th.Mik}. Let us give two such examples to motivate our later axioms.
\bi
\item A straighforward example of random walk on permutations which has large inverted orbits is given by the following asymmetric process driven by rigid permutations: $\tau_n:= s_n \circ \ldots \circ s_1$, with $s_i$ i.i.d. and sampled according to $\mu= \frac 1 d \sum_{i=1}^d \delta_{\sigma_i^{\mathrm{rigid}}}$.  In this case, one clearly has $|O_n|=n+1$ for all $n\geq 0$. This suggests that some kind of isotropy is required to make the question interesting. The natural axiom suggested by Question \ref{q.WZd} is to look at random walks on $\mathfrak{S}(\Z^d)$ with i.i.d increments sampled according to a \textbf{symmetric} probability measure $\mu$ on $\mathfrak{S}(\Z^d)$ (equipped with the  metrizable topology $\calT$ from Definition \ref{d.wobbling}). By symmetric we mean here that for any Borel set $A \subset \mathfrak{S}(\Z^d)$, we require $\mu(A)=\mu(A^{-1}:=\{ \sigma, \sigma^{-1}\in A \})$. This symmetry condition is the same as the symmetry condition required in Kesten's criterion. 
\item Also, if one allows i.i.d increments with \textbf{long-range behavior} (i.e. permutations $\sigma$ such that $\alpha_R(\sigma) = \sup_{\|x\|_2 \leq R}\{ \|\sigma(x) - x\|_2\}$ grows fast with $R$), then it is not hard to come up with symmetric counter-examples. Here is one possible way to build a long-range random walks on $\mathfrak{S}(\Z^d)$: consider any arbitrary injective map from a 3-regular tree $\T^3$ into $\Z$ ($d=1$ is enough here). Now let $\mu$ be the pushforward measure on $\mathfrak{S}(\Z)$ of the uniform measure $\mu_0$ on \textbf{matchings} of $\T^3$ (see for example \cite{LyonsPeres}). Let $O_n\subset \T^3$ be the inverted orbit of the RW on $\mathfrak{S}(\T^3)$ induced by $\mu_0$. It is a rather interesting warming-up problem to check that the inverted orbits $\{O_n\}_{n\geq 0}$ do not satisfy item $i)$ nor $ii)$. In particular, the pushforward process on $\mathfrak{S}(\Z)$ driven by $\mu$ does not satisfy these items either. 
\ei 

\subsection{Classes of diffusion on $\mathfrak{S}(\Z^d)$ considered}

Given the above discussion, we identify four classes of diffusions on $\mathfrak{S}(\Z^d)$ from the most general ones (A) to the most restritive ones (D). Our main results will deal with (A) and (B) only.

\begin{definition}[Four classes of diffusions on $\mathfrak{S}(\Z^d)$]\label{d.classes}
We shall consider random processes $\tau_n := s_n \circ \ldots \circ s_1$, where $\{s_i\}_{i\geq 1}$ are i.i.d and sampled according to a \textbf{symmetric} and \textbf{possibly diffuse} Borel probability mesure $\mu$ on the space $(\mathfrak{S}(\Z^d), \calT)$.  (See definition \ref{d.wobbling}). 

In particular, the family of processes $\{ (n \mapsto \tau_n(\{x\}))_{n\geq 0}\}_{x\in \Z^d}$ are coupled random walks in common conductances given by 
\begin{align}\label{}
C_\mu(x,y):=\mu(\{\sigma, \sigma(x) = y\}) + \mu(\{\sigma, \sigma(y)= x\})
\end{align}

We will consider the following classes of diffusions on $\mathfrak{S}(\Z^d)$. 
\bnum
\item[\textbf{(A)}] \textbf{$L^1$-diffusions on $\mathfrak{S}(\Z^d)$.}

\ni
There exists $k\in L^1(\Z^d)$, s.t. for all $x,y \in \Z^d, \frac{C_\mu(x,y)}{C_\mu(x)} \leq  k(x-y)$. This hypothesis prevents long-range behaviour as in the second example above.

\item[\textbf{(B)}] \textbf{Bounded-range diffusions on $W(\Z^d)$.} 

\ni
We require $k$ to be of finite support, or equivalently that the probability measure $\mu$ is supported on some $W_r(\Z^d) \subset W(Z^d)$ ($r\in \N_+$). This class will correspond to the later notion of \textbf{diffuse-extensive-amenability.} See Definition \ref{d.DEA}.

\item[\textbf{(C)}] \textbf{Finitely-generated diffusions on $W(\Z^d)$}

\ni
One assumes that $\mu$ is finitely supported on $W(\Z^d)$. This is the same as the notion of 
\textbf{extensive-amenability} from \cite{Mik1,Mik2} (see also question \ref{q.WZd}).

\item[\textbf{(D)}] \textbf{``IET'' diffusions on $W(\Z^d), d\geq 1$} 

\ni
The measure $\mu$ is supported on a finite set $\{\sigma_i^{\pm 1}\}_{1\leq i \leq k}\subset W(\Z^d)$ where the permutations $\sigma_i$ are produced through the above recipe from some $\{g_1,\ldots, g_k\} \subset \mathrm{IET}$.  

\enum

Clearly, one has 
\[
(A) \;\supsetneq \;(B) \; \supsetneq \; (C) \; \supsetneq \; (D)\,,
\]
and a positive answer to Question \ref{q.WZd} for any of these classes would imply the amenability of $\mathrm{IET}$. 
\end{definition}

\subsection{Diffuse-extensive-amenability and main results}

Our first main result is to introduce a broad family of random walks $(\tau_n)_{n\geq 0}$ on $\mathfrak{S}(\Z^d)$ which belong to the class \textbf{(A)} from Definition \ref{d.classes} and which do not satisfy items $i)$ or $ii)$ from Theorem \ref{th.Mik}.

\begin{theorem}[See Theorem \ref{th.main}]\label{}
For any $d\geq 3$, there are symmetric probability measures $\mu$ on the space $(\mathfrak{S}(\Z^d),\calT)$ whose associated conductances $C_\mu(x,y)$ defined in Definition \ref{d.classes} induce transient random walks on $\Z^d$ and satisfy for any $x,y \in \Z^d$,
\begin{align*}\label{}
\frac{C_\mu(x,y)}{C_\mu(x)} \leq B\, e^{- b \|x-y\|_2} \text{ for some $b,B>0$.}
\end{align*}
Furthermore, if $(\tau_n)_{n\geq 0}$ is the random walk on $\mathfrak{S}(\Z^d)$ induced by $\mu$, then its inverted orbits $\{O_n\}_{n\geq 0}$ satisfy the following quantitative bound
\begin{align}\label{e.QB}
e^{- n\log(2) \Pb{T=\infty} +o(n)} \leq \Eb{\left(\frac 1 2 \right)^{|O_n|}} \leq e^{-\frac {n+1} 2 \Pb{T=\infty}}\,,
\end{align}
where $\Pb{T=\infty}$ is the probability that the random walk in the conductances $C_\mu$, $n\mapsto \tau_n(\{0\})$, never returns to the origin.  

When $d=2$, one can also find symmetric probability measures $\mu$ on $(\mathfrak{S}(\Z^2),\calT)$ which belong to the class \textbf{(A) of $L^1$-diffusions} from Definition \ref{d.classes} and which do not satisfy items $i)$ or $ii)$ from Theorem \ref{th.Mik}. See Corollary \ref{c.example}. 
\end{theorem}

These random walks on $\mathfrak{S}(\Z^d)$ will be constructed using a well-known continuous time Markov process $t \in \R_+ \mapsto \tau_t \in \mathfrak{S}(\Z^d)$ called the \textbf{random stirring process}. This càdlàg process in $\mathfrak{S}(\Z^d)$ is a process used to analyze/build the so-called \textbf{symmetric exclusion process} on $\Z^d$, which is one of the most studied random particle systems. We will introduce the random stirring process in Section \ref{s.RS} below. See also \cite{Liggett}. 

\medskip

Our second main result is to modify the above processes in order to obtain random walks in class \textbf{(B)}, that is to say with i.i.d and bounded-range increments (i.e. $\mu$ is supported on some $W_r(\Z^d)$). These will be obtained from the above processes by using an appropriate spatial cut-off.

\begin{theorem}\label{th.WZd}
For any $d\geq 3$ and any $r\geq \ceil{3 \sqrt{d}}$ there are symmetric probability measures $\mu$ on $W_r(\Z^d)$ whose associated conductance $C_\mu$ is transient and which induce random walks $(\tau_n)_{n\geq 0}$ on $W(\Z^d)$ satisfying the same exponential decay as in~\eqref{e.QB}. 
\end{theorem}

Finally we generalize the above analysis for euclidean lattices $\Z^d$ to the case of general bounded degree connected graphs $G=(X,E)$. We already defined the {\bf wobbling groups} $W(X)$ and $W_r(X)$ earlier in Definition \ref{d.wobbling}. As suggested to us by M. de la Salle, we introduce the following notion which strengthens in the present setting the notion of {\bf extensive-amenability} from \cite{Kate,Mik1,Mik2}.  

\begin{definition}\label{d.DEA}
Let $G=(X,E)$ be a connected bounded-degree graph.  We will say that the action $W(X) \action X$ is \textbf{diffuse-extensively amenable} if and only if for any $r\in \N_+$, any $x\in X$ and any (possibly diffuse) Borel symmetric probability measure $\mu$ supported on the compact $W_r(X)$, the induced random walk $(\tau_n)_{n\geq 0}$ on $W(X)$ has inverted orbits $\{O_n(x)\}_{n\geq 0}$ which satisfy item $i)$ or $ii)$ from Theorem \ref{th.Mik}.  (In this setting $O_n(x):= \bigcup_{0\leq j \leq n} \tau_j^{-1}(\{ x\})$). 
\end{definition}

\begin{remark}\label{r.DEA} $ $
\bnum
\item Clearly, if a bounded-degree graph is diffuse-extensively amenable, then it must be extensively amenable (according to the equivalent definitions in Question \ref{q.WZd} above). We do not know any examples for which the reverse direction does not hold. See Question \ref{q.EQ}.
\item More generally, if $G$ is a topological group acting faithfully on a countable set $X$, we may wish to define the action $G\action X$ to be \textbf{diffuse-extensively amenable} iff for any symmetric, possibly diffuse, and compactly supported Borel probability measure  $\mu$ on $G$, the induced random walk on $\mathfrak{S}(X)$ has inverted orbits satisfying items $i)$ or $ii)$ from Theorem \ref{th.Mik}. Note that our above definition in the case where $G=W(\Z^d)$ looks slightly more restrictive as we require $\mu$ to be supported on one of the specific compact sets $W_r(\Z^d)$. The reason for this choice is the dichotomy recurrence v.s. transience in the Theorem below. 
\enum
\end{remark}

We will prove the following characterisation of diffuse-extensive-amenability:
\begin{theorem}\label{th.charact}
Let $G=(X,E)$ be a  bounded degree connected graph. The action $W(X) \action X$ is diffuse-extensively amenable if and only if $X$ is recurrent.

(In particular the Euclidean lattices $\Z^d, d\geq 3$ are not diffuse-extensively amenable).  
\end{theorem}

\begin{remark}\label{}
We will not show the implication ``$X$ recurrent'' implies ``$W(X) \action X$ is diffuse-extensively amenable'' as the proof is exactly the same as the proof in \cite{Mik2} that ``$X$ recurrent'' implies ``$W(X) \action X$ is extensively amenable''.
\end{remark}

\ni
{\em Organization.} 
The paper is organized as follows. In Section \ref{s.RS} we give some relevant background on the {\em random stirring process}. In Section \ref{s.proof}, we prove a large deviation bound (estimate~\eqref{e.QB}) on its inverted orbits of sublinear size. The main tool there will be a kind of FKG inequality due to Liggett. 
In Section \ref{s.cutoff}, we add suitable cut-offs on $\Z^d$ to the classical random stirring process in order to obtain processes with i.i.d. increments in $W_r(\Z^d)$. We then axiomatize this procedure for general bounded degree graphs $X$. Finally we end with a conjecture and some open questions suggested by this work. 

\smallskip

\ni
\textbf{Acknowledgments:}

I wish to warmly thank Mikael de la Salle for introducing me to this problem, for numerous enlightening discussions around the concept of extensive-amenability and also for his careful reading of the manuscript. I would also like to thank the anonymous referee for a careful reading as well as Guillaume Aubrun, Kate Juschenko, Nicolas Monod, G\'abor Pete and Fabio Toninelli  for useful discussions and comments on the manuscript. 
This research is supported by the 
ANR grant \textsc{Liouville} ANR-15-CE40-0013 and the ERC grant LiKo 676999.

%
%
%

\section{The random stirring process and its inverted orbit}\label{s.RS}


%
\subsection{The random stirring process on $\Z^d$}\label{ss.RSP}

 We start by defining briefly this process (see for example \cite{Liggett} for background on this process as well as its celebrated companion, the symmetric exclusion process).
 
Let $\{p(x,y)\}_{x,y}$ be a symmetric transition kernel on $\Z^d$, i.e.
\bi
\item[i)] $p(x,y) =p(y,x)$ for all $x,y \in \Z^d$,
\item[ii)] $p(x,y) \geq 0$, 
\item[iii)] $\sum_y p(x,y)=1$ for all $x\in \Z^d$.  
\ei

\begin{figure}[!htp]
\begin{center}
\includegraphics[width=0.8\textwidth]{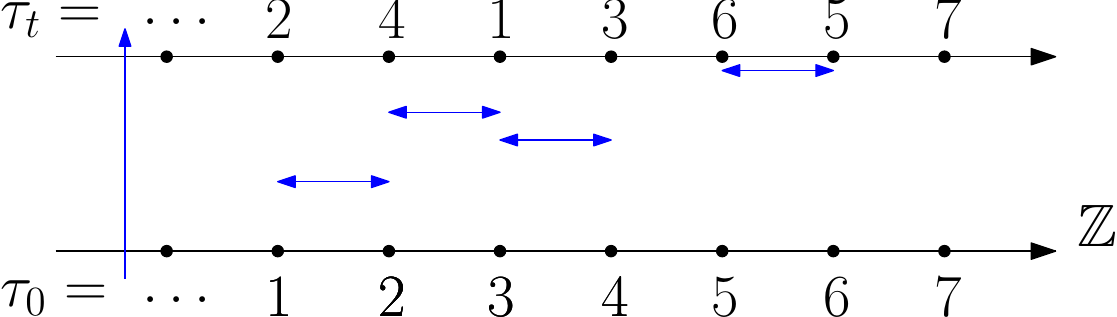}
\end{center}
\caption{(Nearest-neighbour) random stirring process on $\Z$: as time $t$ increases, neighbouring sites are exchanged at rate $\lambda$ (blue arrows). This defines a càdlàg Markov process $t\mapsto \tau_t \in \mathfrak{S}(\Z)$.}\label{f.RSP}
\end{figure}

\begin{definition}\label{d.SP}
Given a symmetric transition kernel $p$ and a real number $\lambda>0$, the {\bf random stirring process} on $\Z^d$ with kernel $\lambda p(\cdot,\cdot)$ is the Markov process $t\mapsto \tau_t$ with values in $\mathfrak{S}(\Z^d)$ which is defined informally as follows: 
\bi
\item $\tau_0 := \mathrm{Id}$ 
\item Independently for each pair $\{x,y\}$, one associates iid exponential clocks of rate $\lambda p(x,y)$ (in other words, independently for each pair $\{x,y\}$, one associates a Poisson Point Process of rate $\lambda p(x,y)$ on $\R_+$). At each time $t$ where an exponential clock rings for a given pair $\{x,y\}$, we define $\tau_t:= T_{xy}\circ \tau_{t-}$ where $T_{xy}$ denotes the transposition of $\{x,y\}$ and $\tau_{t-}$ denotes the left limit of the permutation at the time of the ring. 
\ei
As for any $y\in \Z^d$, we have $\sum_{x\in \Z^d} p(x,y)=1<\infty$, it is a classical fact that one can rigorously make sense of the above construction (see \cite{Liggett}). 

This construction can be defined more generally for any conductance profile $\{c(x,y)\}_{x,y\in \Z^d}$ (recall Definition \ref{d.conductance}) satisfying $\sup_x c(x) = \sup_x \sum_y c(x,y) <\infty$
\footnote{We distinguish the conductances $c(\cdot,\cdot)$ for the continuous time-Markov processes from the conductances $C(\cdot,\cdot)$ in Definition \ref{d.conductance} for discrete time random walks. The relationship between the two reads as follows: $C= e^{L_c}$, where $L_c$ is the symmetric operator defined by $L_c(x,y):=-c(x)(=\sum_y c(x,y))$ if $x=y$ and $L_c(x,y):=c(x,y)$ otherwise.}. 
See figure \ref{f.RSP} for an illustration of the \textbf{random stirring process} on $\Z$. 
\end{definition}

\subsection{Inverted orbits} 

\begin{definition}\label{}
Let $\{p(x,y)\}_{x,y}$ be a symmetric transition kernel on $\Z^d$ and $\lambda>0$ be some fixed positive rate. Let $(\tau_t)_{t\geq 0}$ be the random stirring process induced by $\lambda p(\cdot,\cdot)$.  For each $t\geq 0$, the \textbf{inverted orbit} of the origin $0$ is the random subset $\bar O_t \subset \Z^d$ defined as 
\begin{align*}\label{}
\bar O_t = \bar O_t^\lambda:= \{x\in \Z^d, \exists s\leq t,\,  \tau_s(x)=0\} = \bigcup_{0\leq s \leq t} \tau_s^{-1}(\{0\})\,.
\end{align*}
We also define for any $n\in \N_+$, the inverted orbit along integer times,
\begin{align*}\label{}
O_n = O_n^\lambda:= \{x\in \Z^d, \exists k \in \{0,1,\ldots, n\},  \tau_k(x)=0\} = \bigcup_{0\leq k \leq n} \tau_k^{-1}(\{0\})\,.
\end{align*}
These definitions readily extend to the more general case of a conductance profile $\{c(x,y)\}$. 
\end{definition}

\begin{remark}\label{}
Clearly one always has $O_n \subset \bar O_n$. In particular it is easier for $O_n$ to be of small size than for $\bar O_n$. 
\end{remark}

The motion of one particle, say the origin, under the above random stirring process is simple to describe: it is a continuous time random walk $(\bar X_t^\lambda)_{t\geq 0}$ on $\Z^d$ which jumps at rate $\lambda$ according to the kernel $p(\cdot,\cdot)$. (N.B. in the general case of a conductance profile $\{c(x,y)\}$, it jumps at rate $c(x)$ according to the kernel $\frac{c(x,y)}{c(x)}$). Let $(X_n^\lambda)_{n\in \N}$ be the trace of that continuous time RW discrete time on integer times. In other words, $(X_n^\lambda)_{n\in \N}$ is the discrete time random walk on $\Z^d$ with transition kernel 
\begin{align*}\label{}
[e^{\lambda (p-\mathrm{Id})}](x,y) = e^{-\lambda} \sum_{k\geq 0}  \frac {\lambda^k} {k!} p^{k}(x,y)\,.
\end{align*}
$L:=p-\mathrm{Id}$ is called the \textbf{generator} of this continuous-time random walk. Let us assume that $X_n^\lambda$ starts at the origin and consider the stopping time 
$T:=\inf\{ n \geq 1, X^\lambda_n =0 \}$ which is the first return time to the origin. We shall use throughout the quantity $\FK{\lambda p}{0}{T=\infty}$, which is the probability that $(X_n^\lambda)_{n\geq 0}$ escapes to infinity before returning to the origin. For a general conductance profile $c(\cdot,\cdot)$, we shall use the notation $\FK{c}{0}{T=\infty}$. 

\begin{remark}\label{}
If $\hat X_n$ is the discrete time random walk on $\Z^d$ induced by the transition kernel $p(\cdot,\cdot)$, and if $\hat T:= \inf\{ n\geq 1, \hat X_n =0\}$, then it is easy to check that  
\begin{align*}\label{}
\FK{\lambda p}{0}{T=\infty} \geq (1-e^{-\lambda})\FK{}{0}{\hat T =\infty} = (1-e^{-\lambda}) \frac 1 {R_{eff}(\Z^d, p(\cdot,\cdot))}\,,
\end{align*}
where $R_{eff}(\Z^d, p(\cdot,\cdot))$ is the effective resistance from 0 to infinity of the graph $\Z^d$ with conductances given by $\{p(x,y)\}_{x,y}$.
The above inequality follows from the fact that the continuous-time random walk $(X_n^\lambda)$ will jump during the first unit interval $[0,1)$ with probability $1-e^{-\lambda}$ and will never get back to the origin with probability $\FK{}{0}{\hat T =\infty}$. This lower bound is sub-optimal when $\lambda$ is large: for example in the case of the simple random walk on $\Z^3$, one has $\FK{\lambda p}{0}{T=\infty} \to 1$ as $\lambda \to \infty$ while it is not the case of the R.H.S above. 
\end{remark}

\subsection{Quantitative exponential decay}
\begin{theorem}\label{th.main}
Let $\{p(x,y)\}_{x,y}$ be a symmetric transition kernel on $\Z^d$ and $\lambda>0$ be some fixed positive rate. Then for any $n\geq 0$,
\begin{align}\label{}
\Eb{\left(\frac 1 2 \right)^{|O_n|}} \leq e^{-\frac {n+1} 2 \FK{\lambda p}{0}{T=\infty}}\,.
\end{align}
More generally, for any $p\in[0,1]$ and any conductance profile $\{c(x,y)\}_{x,y}$, we have for all $n\geq 0$, 
\begin{align}\label{e.BH}
\Eb{\left(1-p \right)^{|O_n|}} \leq e^{- p (n+1)  \FK{c}{0}{T=\infty}}\,.
\end{align}
\end{theorem}

\begin{remark}\label{}
As it was pointed out to us by M. de la Salle, this statement is not very far from being optimal: indeed it is not hard to check that $\Eb{|O_n|} = \sum_{k=0}^n \FK{\lambda p}{0}{T>k}$ (see Lemma 4.3 in \cite{Mik2}), 
 which readily implies that 
$\Eb{|O_n|} = n \FK{\lambda p}{0}{T=\infty} + o(n)$. By Jensen, this implies the following lower bound,
\begin{align}\label{}
e^{- n \log(2) \FK{\lambda p}{0}{T=\infty} +o(n)} \leq \Eb{\left(\frac 1 2 \right)^{|O_n|}}\,.
\end{align}
In terms of exponent, we are thus very close as $\log(2)\approx 0.693$ is not far from $1/2$. In fact it is not hard to check that our second estimate ~\eqref{e.BH} is optimal as $p\to 1$. 
Finally, note that the above observation which relates $\Eb{|O_n|}$ with $n\FK{\lambda p}{0}{T=\infty}$ also shows that in the setting of Theorem \ref{th.main}, the inverted orbits are exponentially unlikely to be sublinear if and only if the symmetric random walk $(\hat X_n)_{n\geq 0}$ induced by $p(\cdot,\cdot)$ is transient. 
\end{remark}

The above Theorem has the following straightforward Corollaries.
\begin{corollary}
The following properties on the inverted orbits are satisfied:
\bnum
\item The continuous-time inverted orbit satisfies for any $t>0$, 
\[
\Eb{\left(\frac 1 2 \right)^{|\bar O_t|}} \leq e^{-\frac {\ceil{t}} 2 \FK{\lambda p}{0}{T=\infty}}
\]
\item $\exists a \in(0,1)$, $\Pb{|O_n| < a n} \leq  e^{- a n}$ for all $n\geq 1$. 
\enum
\end{corollary}

\begin{corollary}\label{c.example}
The inverted orbits are exponentially unlikely to be sublinear for the following examples.
\bnum
\item If $d\geq 3$ and $p(x,y):=\frac 1 {2d} 1_{x\sim y}$ on $\Z^d$, exponential decay  holds for any $\lambda>0$. 
\item If $d\in\{1,2\}$ and $p(x,y) \propto \frac 1 {\|x-y\|_2^{d+\alpha}}$, then the exponential decay holds for any $\lambda>0$ if and only if $\alpha \in (0,d)$ . See for example Appendix B.1 in \cite{Caputo}. 
\enum 
\end{corollary}

\begin{remark}\label{}
We highlighted in the above statements the particular case of transition kernels $p(\cdot,\cdot)$ rather than the more general situation of conductance profiles $c(\cdot,\cdot)$ to make these statements more readable for readers unused to continuous-time Markov processes. Still, one may wonder why we added an additional parameter $\lambda>0$. Indeed by construction, our continuous orbits satisfy the obvious scaling identity $|\bar O_t^\lambda| \overset{(law)} = |\bar O_{\lambda t}|$ which makes the parameter $\lambda$ useless from that point of view. There are two reasons. 
\bnum
\item First, the question related to the IET-group is about discrete random walks in $W(\Z^d)$ and in that setting, playing with both $p(\cdot,\cdot)$ and $\lambda$ provides a larger class of examples.  (This is reminiscent of the fact that $C\equiv \lambda C$ for discrete-time walks while $c \not\equiv \lambda c$ for continuous-time RW as it corresponds to a time-change).
\item Second, even though it is in principle harder to obtain large deviation bounds on $O_n \subset \bar O_{t=n}$, our proof in Section \ref{s.proof} will require to focus extensively on the discrete orbits. 
\enum
\end{remark}

\section{Proof of the large deviation bound on the inverted orbits of exclusion process}\label{s.proof}

\subsection{Liggett's inequality in inhomogeneous time}

The main tool in our proof will be the following result due to Liggett which we state in a slightly generalised form (i.e. in inhomogeneous-time).
\begin{proposition}[Proposition 1.7, chapter VIII in \cite{Liggett}]\label{th.L}
Let $L_1,L_2,\ldots, L_k$ be the generators (see \cite{Liggett}) of the random stirring processes on a graph $G=(V,E)$ with respective symmetric transition kernels $p_i: V\times V \to[0,1]$.
Let $\lambda_1, \lambda_2,\ldots,\lambda_k$ be positive rates and let $t_0=0 < t_1<\ldots <t_k=t$.
We consider the random stirring process which evolves at rate $\lambda_1$ according to $L_1$ on the time interval $[0,t_1)$ and then at rate $\lambda_2$ according to $L_2$ on $[t_1, t_2)$ and so on. For any $n\geq 1$, we consider the following diffusion operator on functions $f : V^n \to \R$, 
\begin{align*}\label{}
V_n(t) f := e^{t_1 \lambda_1 L_1} e^{(t_2-t_1) \lambda_2 L_2} \ldots e^{(t-t_{k-1}) \lambda_k L_k}\, f\,.
\end{align*}
I.e. $V_n(t) f(x_1,\ldots,x_n):= \EFK{}{\vec x}{f(\tau_t(x_1), \tau_t(x_2), \ldots, \tau_t(x_n))}$, where $(\tau_t)_{t\geq 0}$ is the (inhomogeneous in time) random stirring process on $G=(V,E)$ induced by $\lambda_1 L_1, \ldots \lambda_k L_k$ and $t_1,\ldots,t_k$ and where
\[
\vec x =(x_1,\ldots, x_n) \in T_n:= V^n \setminus \{\vec y, \exists i \neq j, y_i = y_j \}\,.
\]
Then, we have for any $n\geq 1$ and any bounded \textbf{positive semidefinite}\footnote{In this paper, following \cite{Liggett}, a bounded symmetric function $f:V^n \to \R$ will be called \textbf{positive semidefinite} if it is a positive semidefinite function of each pair of variables. I.e. for any $1\leq i<j\leq n$, any $(x_1,\ldots,x_n)\in V^n$ and any $\{\beta(x)\}_{x\in V}$ with finite support, one has $\sum_{x,y\in V} \beta(x)\beta(y) f(x_1,\ldots,x_{i-1},x, \ldots,y,x_{j+1},\ldots,x_n) \geq 0$.} symmetric function $f:V^n \to \R$,
\begin{align*}\label{}
V_n(t) f(\vec x) \leq U_n(t) f(\vec x) \;\; \text{ for all } \vec x \in T_n\,,
\end{align*}
where $U_n(t)$ is the (inhomogeneous-time) heat semi-group on $G=(V,E)$ induced by independent random walks with same marginal law, i.e. 
\begin{align*}\label{}
U_n(t) f  = [\bigl(e^{t_1 \lambda_1 (p_1-\mathrm{Id})} e^{(t_2-t_1) \lambda_2 (p_2-\mathrm{Id})} \ldots e^{(t-t_{k-1}) \lambda_k (p_k-\mathrm{Id})}\bigr)^{\otimes n}] f\,.
\end{align*}
\end{proposition}

\ni
{\em Proof.}
We refer to \cite{Liggett} for the proof in the time-homogeneous setting. We just point out that in this time-homogeneous case, the hypothesis that the transition kernel $p$ is symmetric is fundamental in the proof. Note here that the random walk corresponding to $e^{t_1 \lambda_1 (p_1-\mathrm{Id})} e^{(t_2-t_1) \lambda_2 (p_2-\mathrm{Id})} \ldots e^{(t-t_{k-1}) \lambda_k (p_k-\mathrm{Id})}$ is no longer symmetric. As such this statement does extend Liggett's result to a considerably larger class of random processes on permutations. Yet, the proof of this extension is immediate: let us argue in the case of two generators, $L_1$ on $[0,t_1)$ and $L_2$ on $[t_1, t_2]$. Then for any bounded positive semidefinite $f : V^n \to \R$, one has $V_n(t) f = V_n^1 (V_n^2 f ) \leq V_n^1 (U_n^2 f)$ since $V_n^1$ preserves the positivity as well as $T_n$. Now (see \cite{Liggett}), it is easy to check that if $f$ is positive semidefinite, so is $U_n^2 f$. In particular, by the homogeneous case, $V_n^1( U_n^2 f) \leq U_n^1(U_n^2 f) = U_n(t) f$ on $T_n$. \qed

\medskip
This very useful inequality has been used  extensively to study properties of the symmetric exclusion process. See for example the use of this inequality in \cite{Arratia} which shares similarities with our present proof (it already relies in fact on the above inhomogeneous version). We will use this inequality as follows:
\begin{corollary}\label{c.L}
Let $A\subset V$ and $f(x_1,\ldots,x_n):= \prod_{i=1}^n 1_{x_i \in A}$. This function is clearly positive semidefinite in each pair of variables. Applying the above Proposition immediately leads to the following inequality: for all $\vec x=(x_1,\ldots, x_n) \in V^n \setminus \{\vec y, \exists i \neq j, y_i = y_j\}$, 
\begin{align*}\label{}
\FK{}{}{\tau_t(x_i) \in A,\text{ for all } 1\leq i \leq n} \leq \prod_{i=1}^n \FK{}{}{Z_t^{x_i} \in A}\,,
\end{align*}
where $\{Z^{z_i}_t\}$ are independent time-inhomegenous walks starting from $x_i$ and driven by the above transition kernels $p_1, \ldots ,p_k$ at rates $\lambda_1,\ldots, \lambda_k$.  
\end{corollary}
%
%

%

\subsection{Bernoulli marking}

The first idea here, which goes back to the proof of (i) in Theorem 4.1 in \cite{Mik2} is to interpret $\Eb{\left(\frac 1 2 \right)^{|O_n|}}$ as the probability that a random subset $A_n \subset \tau_n(O_n)$ is empty where each point in $\tau_n(O_n)$ belongs to $A_n$ independently with probability $1/2$ (indeed $|\tau_n(O_n)|=|O_n|$). There is a very convenient and dynamical way to sample the random subset $A_n \subset \tau_n(O_n)$ which goes as follows (once again suggested by \cite{Mik2}): let $s_1, s_2, \ldots, s_n$ the permutations on $\Z^d$ induced by the (independent) random stirring processes along the intervals $(0,1), (1,2), \ldots (n-1,n)$, let also $t_0,t_1,\ldots,t_n$ be $(n+1)$ independent standard Bernoulli variables in $\{0,1\}$. 
\bnum
\item At time $0$, let $A_0 = \emptyset$ if $t_0=0$ and $A_0:=\{0\}$ otherwise. 
\item At time $1$, let $A_1:= s_1(A_0)\setminus \{0\}$ if $t_1=0$ and $A_1:=s_1(A_0) \cup \{0\}$ otherwise. 
\item We keep going inductively for $i\geq 2$, let $A_i:= s_i(A_{i-1})\setminus \{0\}$ if  $t_i=0$,  and $A_i:=s_1(A_{i-1}) \cup \{0\}$ otherwise.
\enum
As such, clearly $A_n$ is a uniform subset of $\tau_n(O_n)$ and thus $A_n =\emptyset$ with probability exactly $(1/2)^{|O_n|}$.

\subsection{Nearly Bernoulli marking via large $N$-reservoirs}
\begin{figure}[!htp]
\begin{center}
\includegraphics[width=\textwidth]{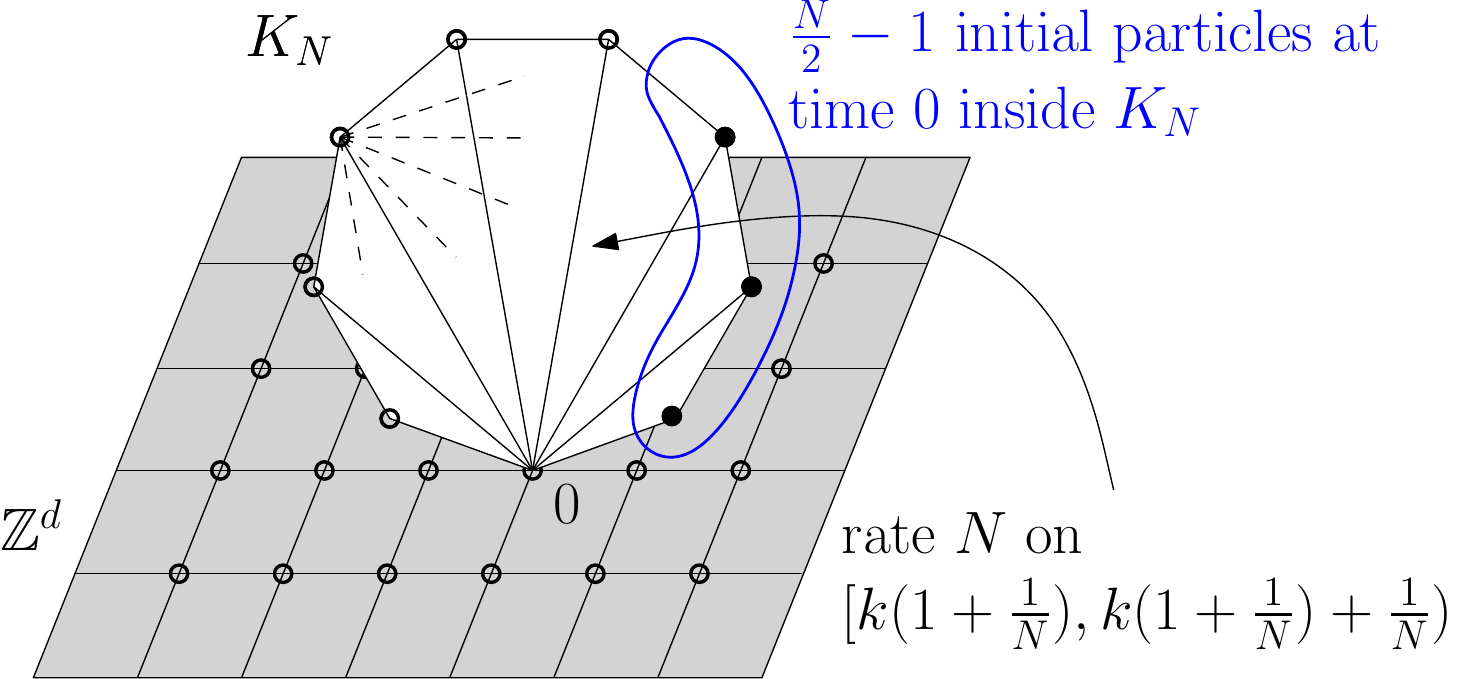}
\end{center}
\caption{}\label{f.KN}
\end{figure}

We will not achieve such a perfect Bernoulli marking, but we shall instead get as close as we want by relying on the following approximation procedure (see Figure \ref{f.KN}). The very important feature of our approximation will be that not only it will serve as a nearly perfect Bernoulli(1/2)-marking but it will also evolve as a time-inhomogeneous exclusion process on a larger graph than $\Z^d$, 
where one connects a large ``reservoir'', namely a complete graph $K_{N}$ with  $N\gg 1$, to the origin $0\in \Z^d$. See Figure \ref{f.KN}. This second property will allow us to reduce the analysis to independent time-inhomogeneous random walks running on $\Z^d \cup K_{N}$ thanks to Proposition \ref{th.L}.

We will work for each even $N\in 2 \N_+$ on the graph $G_N$ with vertex set $\Z^d \cup K_N$ (where $K_N$ is the complete graph on $N$ vertices) and where we connect $\Z^d$ with $K_N$ by adding one edge between each $N$ vertex of $K_N$ and the origin $0$ of the lattice $\Z^d$. Let us now set up the inhomogeneous rates we shall use on this lattice: 
\bi
\item For each $k\geq 1$, let $I_k$ be the unit interval $[k-1+\frac {k} N, k +\frac {k}N)$. On each such interval $I_k$, the random stirring process is driven by $\lambda p(\cdot,\cdot)$ where $p$ is our initial symmetric kernel on $\Z^d$. In other words, each edge $e=(x,y), x,y \in \Z^d$ switch at rate $\lambda p(x,y)$ on $I_k$. Furthermore, the rates inducing the switches between $K_N$ and $\{0\}$ are frozen, i.e. equal to 0, on these intervals. 
\item For each $k\geq 0$, let $J_k$ be the interval $[k+\frac k N, k + \frac {k+1} N)$. 
On each such interval $J_k$, assume each edge from $K_N$ to $\{0\}$ switch at rate $N$, while the rest of the process (on $\Z^d$) is frozen on these intervals. 
 
\item By symmetry, there is no need to describe the rates between vertices in $K_N$, we may assume these rates to be zero if we wish to. 
\ei 

This construction induces as in Definition \ref{d.SP} two natural processes:
\bnum
\item[a)] First a Markov process $t \mapsto \tau^{(N)}_t$ on the permutations of $\Z^d \cup K_N$ which only permutes sites in $\Z^d$ on intervals $I_k, k\geq 1$ while it only permutes sites between $K_N$ and the origin $0$ on intervals $J_k, k\geq 0$.
\item[b)] It also induces a random walk on the permutations of $\Z^d$ which has exactly the law of the original random stirring process. For each $k\geq 1$, let $s_k$ be the permutation of $\Z^d$ induced by the stirring process restricted to the interval $I_k$. As such $\{ s_k\}_{k\geq 1}$ are i.i.d and $\tau_n:=s_{n} \circ \ldots \circ s_1$ has the same law as in Theorem \ref{th.main}.
\enum
We will use the first process $t \mapsto \tau^{(N)}_t$ to induce a nearly Bernoulli-marking: initialise the system with $\frac N 2 -1$ particles in $K_N$, say at sites $\{x_1,\ldots,x_{\frac N 2 -1}\}$ (and no particles elsewhere). Now for each $t\geq 0$, let 
\begin{align*}\label{}
A_t^{(N)}:= \{ \tau^{(N)}_t(x_i), 1\leq i \leq \frac N 2 - 1\} \cap \Z^d
\end{align*}
This set is a random subset of $\Z^d$ evolving in time. For each $n\in \N$, let $\hat n:= n(1+\frac 1 N)+\frac 1 N$. Our approximation has been constructed so that $A_{\hat n}^{(N)}$ is very close to the exact Benoulli-marked random set $A_n \subset \tau_n(O_n)$ we considered above. The slight technical difficulty here is that given the set $\tau_n(O_n)$, the events $\{x\in A_{\hat n}^{(N)}\}, x\in \tau_n(O_n)$ are no longer independent (in other words $A_{\hat n}^{(N)}$ is not quite a uniform random subset of $\tau_n(O_n)$). Yet, to overcome this lack of independence, the following stochastic domination type of argument allows us to conclude.

\subsection{Stochastic domination}


We will show the following lower bound. 
\begin{lemma}\label{l.LB}
If $N$ is large enough,
\[
\Pb{A_{\hat n}^{(N)} = \emptyset} \geq \Eb{(1/2)^{|O_n|}}
\]
\end{lemma}

\ni
{\em Proof.}
The randomness used to sample the random set $A_{\hat n}^{(N)}$ can be naturally decomposed into two sources of randomness. First, recall $s_1,\ldots, s_{n}$ are the permutations of $\Z^d$ induced by the stirring processes on each of the intervals $(I_i)_{1\leq i \leq n}$. Let us call $\F$ the filtration generated by these permutations. The other source of randomness needed to run $A_{\hat n}^{(N)}$ is the shuffling within $K_N \cup \{0\}$ which arises on each intervals $J_k, k\geq 0$. For each $k\geq 0$, call $\calG_k$ the filtration generated by that shuffling on $J_k$ (i.e. by the stirring process in $K_{N}\cup \{0\}$ restricted to $J_k$).  With these definitions, the permutation $\tau_{\hat n}^{(N)}$ of $\Z^d\cup K_N$ is measurable w.r.t to $\sigma(\F, (\calG_k)_{0\leq k \leq n})$ while the permutation $\tau_n$ of $\Z^d$ is measurable w.r.t $\F$ only. In particular the sets $O_n$ and $\tau_n(O_n)$ are also measurable w.r.t $\F$. We may thus write, 

\begin{align*}\label{}
\Pb{A_{\hat n}^{(N)} = \emptyset} & = \Eb{\Pb{A_{\hat n}^{(N)} = \emptyset \md \F}}  \\
& = \Eb{ \Eb{ \prod_{x\in \tau_n(O_n)} 1_{x\notin A_{\hat n}^{(N)}} \md \F} }
\end{align*}

Now, let us use any prescribed ordering of the set $\tau_n(O_n)$ so that may write $O_n = \{z_1,\ldots, z_m\}$, where $m=|O_n|$. For each $1\leq i \leq m$, let us define $k_i:= \sup \{ j\leq n, \tau_{j}(z_i) = 0 \}$, i.e the last passage in $0$ of the orbit of $z_i$ before time $n$. This random time $k_i\in \{0,\ldots, n\}$ is also measurable w.r.t $\F$ which allows us to go one step further:
\begin{align*}\label{}
\Pb{A_{\hat n}^{(N)} = \emptyset}  & = \Eb{ \Eb{ \prod_{i=1}^{m=|O_n|} 1_{z_i\notin A_{\hat n}^{(N)}} \md \F} } \\
& = \Eb{\Eb{   \Eb{  \prod_{i=1}^{m}1_{z_i\notin A_{\hat n}^{(N)}} \md \F, (\calG_k)_{k\neq k_1}}  \md \F }} \\
& = \Eb{\Eb{ \prod_{i=2}^{m}1_{z_i\notin A_{\hat n}^{(N)}}   \Eb{  1_{z_1\notin A_{\hat n}^{(N)}}  \md \F, (\calG_k)_{k\neq k_1}}  \md \F }}\,,
\end{align*}
since it can be easily checked that the events $\{ z_i \notin A_{\hat n}^{(N)}\}_{2\leq i \leq m}$ are measurable w.r.t $\sigma(\F, (\calG_k)_{k\neq k_1})$. 
We will prove below that 
\begin{align}\label{e.Bel}
\Eb{  1_{z_1\notin A_{\hat n}^{(N)}}  \md \F, (\calG_k)_{k\neq k_1}}  \geq \frac 1 2
\end{align}
which will imply 
\begin{align*}\label{}
\Pb{A_{\hat n}^{(N)} = \emptyset}  & \geq \Eb{\Eb{ \frac 1 2  \prod_{i=2}^{m}1_{z_i\notin A_{\hat n}^{(N)}} \md \F }}\,,
\end{align*}
and by repeating inductively the same argument from $i=2$ to $i=m$, this gives us the desired result, i.e.
\begin{align*}\label{}
\Pb{A_{\hat n}^{(N)} = \emptyset}  & \geq \Eb{\Eb{ (\frac 1 2)^{m} \md \F }} = \Eb{(\frac 1 2 )^{|O_n|}}\,.
\end{align*}
Let us then show the estimate~\eqref{e.Bel}. Given $\F, (\calG_k)_{k\neq k_1}$, $z_1$ will not belong to $A_{\hat n}^{(N)}$ if and only if the shuffling inside $K_n\cup\{0\}$ on the interval $J_{k_1}$ will not end up with a particle at the origin. The probability that there is no particle is larger than the probability $p_{\mathrm{update}}$ that there is at least one update of the origin on $J_k$, times the probability $p_0$ that the last transposition did not place a particle at $\{0\}$. From the construction of our process, $p_{\mathrm{update}}=1-e^{-N \times N /N}=1-e^{-N}$. For the estimate of $p_0$, note that at the time of the last switch, there are at most $\frac N 2-1$ particles in $K_N$, hence $p_0 \geq (\frac N 2 + 1)/N= \frac 1 2 + 1/N$. All together this gives us a lower bound of $(1-e^{-N})(\frac 1 2 +1/N)$ which is higher than $1/2$ for $N$ large enough. This ends our proof. 
\qed

\subsection{Upper bound via independent random walks}
We shall now rely on Corollary \ref{c.L} to show the following domination. 
\begin{lemma}\label{}
\[
\limsup_{N\to \infty} \Pb{A_{\hat n}^{(N)} =\emptyset} \leq e^{-\frac {n+1} 2 \FK{\lambda p}{0}{T=\infty}}\,.
\]
\end{lemma}
\ni
Clearly this Lemma together with Lemma \ref{l.LB} concludes the proof of Theorem \ref{th.main}. 

\ni
{\em Proof.}
Let $\vec x:= (x_1,\ldots, x_{\frac N 2 -1})$ be our initial set of particles in $K_N$ (see Figure \ref{f.KN}). Corollary \ref{c.L} implies that 
\begin{align}\label{e.UB}
\Pb{A_{\hat n}^{(N)}= \emptyset} &= \Pb{\tau_{\hat n}(x_1)\in K_N, \text{ for all }1\leq i \leq \frac N 2 -1} \nn  \\
& \leq \prod_{i=1}^{\frac N 2 -1} \FK{}{x_i}{Z_{t=\hat n} \in K_N} \nn \\
&= \exp\bigl( (\frac N 2 -1) \log \FK{}{x_1}{Z_{t=\hat n} \in K_N}\bigr) \text{  by symmetry} 
\end{align}
We thus need an upper bound on $\FK{}{x_1}{Z_{t=\hat n} \in K_N}$, i.e. a lower bound on the probability $p_{\mathrm{escape}}$ that the particle initially at $x_1$ escapes from $K_N$ at time $t=\hat n$. One possible way of escaping is as follows:
\bnum
\item First, $Z_t^{x_1}$ needs to be placed in $\{0\}$ at the end of one of the shuffling on intervals $J_k$. At the beginning of each of these $J_k$, all the $N$ sites of $K_N$ have the same probability $p_{\mathrm{marked}}$ to be placed at $\{0\}$ at the end of $J_k$ (note that $p_{\mathrm{marked}}$ does not depend on what happened before as it also takes into account the event that a site in $K_N$ without a particle is sent to $\{0\}$). We use below the fact that the motion of one i.i.d particle $t\mapsto Z_t^{x_1}$ has the same law as $t\mapsto \tau_t(x_1)$. In particular, we see that $N\, p_{\mathrm{marked}}$ is exactly the probability that one of the sites initially in $K_N$ is placed at $\{0\}$ at the end of $J_k$ under the exclusion dynamics on $J_k$. Now it is easy to see that this latter probability is greater than the probability that there is at least one update along $J_k$ (i.e. $= 1-e^{-N}$) times the probability that the site initially at $\{0\}$ (at the beginning of $J_k$) is not sent back to $\{0\}$ at the end of $J_k$. By looking at what happens at the last update, this is greater than $1-\frac 1 N$. This gives us $p_{\mathrm{marked}} \geq \frac 1 N(1-e^{-N})(1-\frac 1 N) \geq \frac 1 N - \frac 2 {N^2}$ for $N$ large enough. Now, the probability that $Z_t^{x_1}$ gets placed in $\{0\}$ at the end of one of the intervals $(J_k)_{0\leq k \leq n}$ is exactly
\begin{align*}\label{}
1-(1-p_{\mathrm{marked}})^{n+1} \geq 1-(1-\frac 1 N + \frac 2 {N^2})^{n+1} = \frac{n+1}{N} + O(N^{-2})
\end{align*}
\item Now, if the particle initially at $x_1$ is placed in $\{0\}$ at the end of a $J_k$, one way to get a lower bound is to ask this particle to never get back to the origin. This probability depends on what is $0\leq k \leq n$ but in any case is larger than $\FK{\lambda p}{0}{T=\infty}$
\enum 
All together, we find that 
\begin{align*}\label{}
p_{\mathrm{escape}} \geq \frac{n+1}{N}\FK{\lambda p}{0}{T=\infty} + O(N^{-2})
\end{align*}
Getting back to ~\eqref{e.UB}, this gives us 
\begin{align*}\label{}
 \Pb{A_{\hat n}^{(N)}} & \leq \exp\bigl((\frac N 2 -1) \log(1-\frac{n+1}{N}\FK{\lambda p}{0}{T=\infty} + O(N^{-2}) )\bigr) \\
 & = \exp\bigl(- (\frac N 2 -1)\frac{n+1}{N}\FK{\lambda p}{0}{T=\infty} + O(N^{-1}) )\bigr)\,.
\end{align*}
Now, sending $N\to \infty$ concludes our proof. 
\qed

Note that the exact same proof leads to the more general estimate~\eqref{e.BH} in our main Theorem \ref{th.main}. Indeed, the only change is to consider $\floor{pN}-1$ initial particles in the ``reservoir'' $K_N$ instead of $\frac N 2 -1$. Similarly the same proof applies if one replaces $\lambda p(\cdot,\cdot)$ by any general conductance profile $c(\cdot,\cdot)$.

\section{spatial cut-off and a characterization of diffuse-extensive-amenability}\label{s.cutoff}

\subsection{The case of $\Z^d$}
We start by explaining how one can obtain bounded range increments on $\Z^d$ (which will prove Theorem \ref{th.WZd}).  Clearly, the classical random stirring process $t\mapsto \tau_t$ (i.e. with transition kernel $p(x,y):= \frac 1 {2d} 1_{x\sim y}$) is not of finite range. In other words, for any $t>0,\, \tau_t \notin W(\Z^d)$ a.s. (the process $t\mapsto \tau_t$ belongs to the class \textbf{(A)} of $L^1$-diffusions from Definition \ref{d.classes} but not to the class \textbf{(B)}). To keep the technology used in Section \ref{s.proof} valid (namely, Liggett's inequality), we will still rely on some underlying Exclusion process, but to make increments of bounded range, we will remove some edges in order to confine the dynamics. On the other hand, if we confine things too much then the process will stay in some $W_{r_0}(\Z^d)$ for ever and the inverted orbits will remain very small ($|O_n| \leq C\, r_0^d, \forall n\geq 1$). 

\begin{figure}[!htp]
\begin{center}
\includegraphics[width=0.3\textwidth]{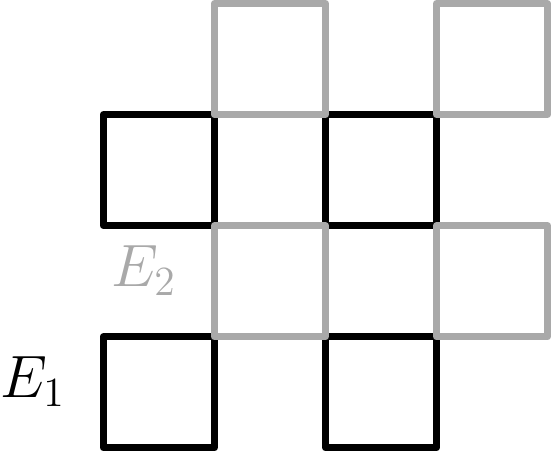}
\end{center}
\caption{Spatial cut-off used on $\Z^2$, i.e. $\E^2=E_1\cup E_2$}\label{f.E1E2}
\end{figure}
We will proceed as follows. First we shall divide the (nearest-neighbour) edges of $\Z^d$ (denoted by $\E^d$) into two sets $E_1 \cup E_2$. Consider for this all the unit cubes $[0,1]^2$ translated by $(2\Z)^d$, i.e. all cubes of the form $(2n_1,\ldots, 2n_d)+[0,1]^d$. Let $E_1$ be the set of all edges belonging to one of these cubes. See Figure \ref{f.E1E2}. And Let $E_2:= \E^d \setminus E_1$. (i.e. $E_2$ is the set of edges adjacent to all $(2\Z)^d$-translates of the cube $[1,2]^d$).
Now we set up the following transition kernels:  let $p_1(x,y):= \frac 1 d 1_{x\sim_{E_1} y}$ and $p_2(x,y):= \frac 1 d 1_{x\sim_{E_2} y}$. 

We introduce the following inhomogeneous-time random stirring process $t\geq 0\mapsto \tau_t \in W(\Z^d)$ : for any $k\in \N$, on the time-interval $[k,k+1/3)$, we let the process evolve according to a time-homogeneous random stirring process with transition kernel $p_1(\cdot,\cdot)$. Then, on $[k+1/3,k+2/3)$, we let the process evolve according to a time-homogeneous stirring process with kernel $p_2$ and finally on $[k+2/3,k+1)$, we let things evolve again according to   the kernel $p_1$. One may write this in terms of operators as follows: if $L_1,L_2$ are the generators of the random stirring processes induced by $p_1,p_2$, then the Markov operator describing the law of $\tau_{t=1}$ is 
\[
K= e^{1/3 L_1} e^{1/3 L_2} e^{1/3 L_1}\,.  
\]

The reason why we split each unit interval $[k,k+1)$ into three intervals instead of $p_1$ on $[k,k+1/2)$ and $p_2$ on $[k+1/2,k+1)$ is for the following reason: for any $x\in \Z^d$, the discrete-time process
\[
n\in \N \mapsto \tau_{n}(\{x\})\,, 
\]
is a RW with transition kernel given by 
\[
P = e^{\frac 1 3 (p_1 - \mathrm{Id})}e^{\frac 1 3 (p_2 - \mathrm{Id})}e^{\frac 1 3 (p_1 - \mathrm{Id})}\,.
\]
In particular, with our choice of splitting, $P$ is symmetric ($P=P^T$) and $n\mapsto  \tau_{n}(\{x\})$ is thus a RW in conductances $C(x,y):=P(x,y)$ (see Definition \ref{d.conductance}). It is straighforward to check from our construction that these conductances satisfy item $e)$ from Proposition \ref{pr.conductance}. This shows that if $d\geq 3$, $\FK{P}{0}{T=\infty}>0$ and the exact same proof as in Section \ref{s.proof} implies that the inverted orbits of the (time-inhomogeneous) dynamics $t\mapsto \tau_t$ satisfy at any integer times $n\geq 0$, 
\[
\Eb{\left(\frac 1 2 \right)^{|O_n|}} \leq e^{-\frac {n+1} 2 \FK{P}{0}{T=\infty}}\,.
\]
(In the present case the above symmetrisation step is not fully necessary but it will be in the more general case below).  
Finally, it is immediate from our construction that $\tau_n=s_n\circ \ldots s_1$ where the random permutations are i.i.d (sampled according to $\tau_{t=1}$)  and belong to $W_{r}(\Z^d)$ with $r=\ceil{3 \sqrt{d}}$.  This ends our proof. (N.B. the probability measure $\mu$ on $W_r(\Z^d)$ which induces this process can be described using the above Markov operator $K$ by $\mu(A):= P_{t=1}(1_A)(\delta_{\mathrm{Id}}) = [K \, 1_{A}](\delta_{\mathrm{Id}}) $).

\subsection{The general case of connected-bounded-degree graphs}
In this Section, we will prove  Theorem \ref{th.charact} which provides a characterisation of our notion of {\em diffuse-extensive-amenability} from Definition \ref{d.DEA}.

\smallskip

\ni
{\em Proof of Theorem \ref{th.charact}.}

\ni
Let $G=(X,E)$ be a connected bounded degree graph. Let $h\in \N_+$ be such that $\mathrm{deg}(x) \leq h\,, \forall x\in X$.  Let $E:=\{e_1,e_2,\ldots, e_n,\ldots \}$ by any ordering of the set of edges. We define recursively the following subsets of $E$:
\bnum
\item First we defined $E_1\subset E$ as follows. Let $E_{1,0}:=\emptyset$. For any $n\geq 0$, let $E_{1,n+1}$ be $E_{1,n} \cup \{e_{n+1}\}$ iff the edge $e_{n+1}$ does not share any common vertex with the edges already in $E_{1,n}$. Otherwise let $E_{1,n+1}=E_{1,n}$.  Define $E_1:= \bigcup E_{1,n}$. 
\item Assume $E_1,\ldots, E_k \subset E$ have been defined. Construct $E_{k+1}$ recursively as follows. Let $E_{k+1,0}:=\emptyset$. For any $n\geq 0$, let $E_{k+1,n+1}$ be $E_{k+1,n} \cup \{e_{n+1}\}$ iff the edge $e_{n+1}$ does not belong to  $\bigcup_{m=1}^k E_m$ and does not share any common vertex with the edges already in  $E_{k+1,n}$. Otherwise let $E_{k+1,n+1}=E_{k+1,n}$.  Define $E_{k+1}:= \bigcup E_{k+1,n}$. 
\item Iterate until one has $\bigcup_{m=1}^k E_m = E$. 
\enum
Let us argue that this procedure ends in finite time (i.e. that there is $k\geq 1$, such that $\bigcup_{m=1}^k E_m = E$). For any $j\geq 1$, let us see what happens for the $j$-th edge $e_j\in E$. Consider the finite set  $F=F_j\subset E$ of all edges  which share a common vertex with $e_j$ and whose labels are $<j$. In particular, we have $|F| \leq 2(h-1)$. Along the process of defining the set $E_1$, if $F\cap E_{1,j-1} =\emptyset$, then by construction $e_j$ will belong to $E_1$. Otherwise, this means at least one egde in $F$ will belong to $E_1$. If we are in the second case, we proceed to the construction of $E_2$. Along that construction, either we add $e_j$ to $E_2$ or this means at least one new edge in $F$ (which did not appear in $E_1$) has been added in $E_2$. As $F$ is finite, this readily implies that the edge $e_j$ needs to be added to $E_k$ at the latest when $k=2(h-1) +1$. (I.e. $E=E_1 \cup  \ldots E_{2h-1}$). 

Now, for any $i\in \{1,\ldots, 2h-1\}$, define the conductance profile $c_i$ on $G=(X,E)$ as follows: for any $x,y\in E$, if the edge $e=\{x,y\} \in E_i$ then define $c_i(x,y):=1$ and otherwise let $c_i(x,y):=0$. (With this definition, there may be some sites $x$ here where $c_i(x)=\sum_y c_i(x,y)=0$ but it is harmless for continuous-time dynamics). 
Similarly as in the above Section, if $L_1,\ldots, L_{2h-1}$ denote the generators of the random stirring processes associated to these conductance profiles, we shall consider the inhomogeneous-time random stirring process which evolves on each unit time interval as follows. Let $N:=2h-1$. For any $k\geq 0$ and any $1\leq i \leq N$, let the process $t\mapsto \tau_t$ evolve on $[k+\frac{i-1}{2N}, k+ \frac {i}{2N})$ according to the random stirring process induced by $L_i$. This defines our stochastic process on all intervals of the form $[k,k+1/2)$. On the remaining intervals, symmetrise exactly as in the above argument in order to induce symmetric random walks $n\mapsto \tau_{n}(x),\,  \forall x\in X$.   
I.e. the Markov operator describing the law of $\tau_{t=1}$ is 
\begin{align*}\label{}
K:=e^{\frac 1 {2N} L_1} \ldots e^{\frac 1 {2N} L_N} e^{\frac 1 {2N} L_N}  \ldots e^{\frac 1 {2N} L_1}\,.
\end{align*}
It is straightforward to check that 
\bi
\item $\tau_{t=1}\in W_{4h-2}(X)$
\item The conductance profile $C$ associated to $K$ satisfies $C(x,y)\geq \delta 1_{x\sim_G y}$ for some $\delta=\delta(h)>0$. 
\ei
As item $e)$ from Proposition \ref{pr.conductance} is also valid for transient bounded-degree graphs, this concludes our proof. 
\qed

%


\section{Conclusion and some open questions}

We have shown that for rather generic examples of diffusions on $W(\Z^d), d\geq 3$, inverted orbits $\{O_n\}_{n\geq 0}$ are exponentially unlikely to be sublinear. Based on this, it would sound very counter-intuitive (and very remarquable!) to us if for ANY symmetric and finitely supported probability measure $\mu$ on $W(\Z^d)$, their corresponding inverted orbits had better-than-exponential probability to be sublinear. We do not believe this can hold and state the following conjecture.

\begin{conjecture}
When $d\geq 3$, the action $W(\Z^d) \action \Z^d$ is not extensively amenable. 
\end{conjecture}

More generally, one may ask the following question.
\begin{question}\label{q.EQ}
Is it the case that for any  bounded-degree connected graph $G=(X,E)$, the action $W(X)\action X$ is extensively amenable if and only if it is diffuse-extensively amenable ?
(If so, by Theorem \ref{th.charact}, $W(X)\action X$ would be extensively amenable if and only if $X$ is recurrent).  
\end{question}

\bibliographystyle{alpha}

\begin{thebibliography}{GG}



\bibitem{Amir}
Amir, G. and Virág, B. 
\newblock Speed exponents of random walks on groups. 
\newblock {\em International Mathematics Research Notices,} 2017(9), 2567--2598. 2017. 


\bibitem{Arratia}
Arratia, R., 
\newblock Symmetric exclusion processes: a comparison inequality and a large deviation result. 
\newblock {\em The Annals of Probability}, pp.53--61. 1985. 


\bibitem{Bartho1}
Bartholdi, L. and Erschler A.
\newblock  Growth of permutational extensions.
\newblock {\em Inventiones mathematicae} 189.2 pp. 431--455, 2012.


\bibitem{Bartho2}
Bartholdi, L. and Erschler A.
\newblock Poisson?Furstenberg boundary and growth of groups. 
\newblock {\em Probability Theory and Related Fields,} 168(1-2), 347--372. 2017.



\bibitem{Caputo}
Caputo, P., Faggionato, A. and Gaudillière, A., 
\newblock Recurrence and transience for long-range reversible random walks on a random point process. 
\newblock {\em Electronic Journal of Probability}, 14, pp. 2580--2616. 2009.

\bibitem{Cornulier}
De Cornulier, Y.  Groupes pleins-topologiques (d?après Matui, Juschenko, Monod,...). {\em Astérisque,} 361, pp.183-223. 2014.


\bibitem{Dahmani}
Dahmani, F., Fujiwara, K. and Guirardel, V., Free groups of interval exchange transformations are rare. {\em Groups, Geometry, and Dynamics,} 7(4), pp.883--911. 2013.

\bibitem{Van}
Van Douwen, E.K. Measures invariant under actions of F2. {\em Topology and its Applications,} 34(1), pp.53--68. 1990


\bibitem{Kate}
Juschenko, K. and Monod, N.
\newblock Cantor systems, piecewise translations and simple amenable groups. 
{\em Annals of mathematics,} 178(2), pp.775--787, 2013.


\bibitem{wobbling}
Juschenko, K. and de la Salle, M., 
Invariant means for the wobbling group.
{\em Bulletin of the Belgian Mathematical Society-Simon Stevin}, 22(2), pp.281--290.
2015.	
	



\bibitem{Mik1} Juschenko, K., Nekrashevych, V. and de la Salle, M.
\newblock Extensions of amenable groups by recurrent groupoids.
\newblock {\em Inventiones mathematicae,} 206(3), pp. 837--867. 2016.


\bibitem{Mik2}
Juschenko, K., Matte Bon, N., Monod, N. and de la Salle, M.
\newblock Extensive amenability and an application to interval exchanges. 
\newblock {\em Ergodic Theory and Dynamical Systems,} pp.1--25, 2016.



\bibitem{Virag}
Kotowski, M. and Virág, B., Non-Liouville groups with return probability exponent at most 1/2. {\em Electronic Communications in Probability,} \textbf{20.} 2015. 



\bibitem{Kumagai}
Kumagai, T.
\newblock Random walks on disordered media and their scaling limits. {\em St-Flour lecture notes.} Springer.  2014.

\bibitem{Liggett}
Liggett, T.M. Interacting particle systems (Vol. 276). {\em Springer.} 2012. 

\bibitem{LyonsPeres}
Lyons, R. and Peres, Y., 
\newblock Probability on trees and networks (Vol. 42). {\em Cambridge University Press.}
2016. 

\bibitem{Paterson}
Paterson, A.L.  Amenability. {\em American Mathematical Soc.} \textbf{(No. 29).} 2000.

\bibitem{Saloff}
Saloff-Coste, L. and Zheng, T., Isoperimetric profiles and random walks on some permutation wreath products. {\em preprint.} arXiv:1510.08830. 2015. 


\bibitem{Viana}
Viana, M.,  Ergodic theory of interval exchange maps. {\em Revista Matemática Complutense,} 19(1), pp.7--100. 2006.




\end{thebibliography}

\end{document}